\documentstyle[11pt]{article}
\renewcommand{\baselinestretch}{1.3}

\newcommand{\labeq}[1]{  \be \label{#1}  }
\newcommand{\labea}[1]{  \bea \label{#1}  }
\newcommand{\ben}{\begin{enumerate}}
\newcommand{\een}{\end{enumerate}}

\newcommand{\bea}{\begin{eqnarray}}
\newcommand{\ba}{\begin{array}}
\newcommand{\bean}{\begin{eqnarray*}}
\newcommand{\ea}{\end{array}}
\newcommand{\eea}{\end{eqnarray}}
\newcommand{\eean}{\end{eqnarray*}}
\newcommand{\be}{\begin{equation}}
\newcommand{\ee}{\end{equation}}

\alph{enumii}
\roman{enumiii}
\newtheorem{thm}{Theorem}[section]
\newtheorem{prop}[thm]{Proposition}
\newtheorem {lem}[thm]{Lemma}

\newtheorem {dfn}{Definition}

\def \lequiv 
{ \;\; \raise-2dd \hbox{$\stackrel{ \prec}{\textstyle \frown}$}\;\;}
\def \gequiv 
{ \;\; \raise-2dd \hbox{$\stackrel{ \succ}{\textstyle \frown}$}\;\;}
\def \Cap {\mbox{\rm Cap}}
\def \E {{\bf E}}
\def \Var {{\bf Var}}
\def \P {{\bf P}}
\def \R {{\bf R}}
\def \Z {{\bf Z}}
\def \dn{{\cal D}_n}
\def \d{\mbox{\boldmath $d$}}
\def\one{\mbox{\bf 1}}
\def\Cox{\hfill \Box}
\def\dist{{\rm dist}}

\title{The Trace of Spatial Brownian Motion is Capacity-equivalent to the 
Unit Square}
\author{Robin Pemantle$^1$
\and Yuval Peres$^{2,3}$
 \and Jonathan W. Shapiro$^{2,4}$ } 
\begin{document}
\begin{titlepage}
\maketitle
\addtocounter{footnote}{1}
\footnotetext{Department 
of Mathematics, University of Wisconsin, Madison, WI 53706.
  Supported in part by 
National Science Foundation grant \# DMS 9300191, by a Sloan Foundation
Fellowship, and by a Presidential Faculty Fellowship.}
\addtocounter{footnote}{1}
\footnotetext{Department of Statistics, University of
 California, Berkeley, California 94720.}
\addtocounter{footnote}{1}
\footnotetext{Research partially supported by NSF grant 
\# DMS-9404391.}
\addtocounter{footnote}{1}
\footnotetext{Research supported by a Line and Michel  Lo\`{e}ve
Fellowship.}

\begin{abstract}
  We show that
 with probability 1, the trace $B[0,1]$ of Brownian motion
 in space,
has positive capacity with respect
to exactly the same kernels as the unit square.
More precisely, the energy 
of occupation measure 
on $B[0,1]$  in the kernel $f(|x-y|)$, is
bounded above and below
by constant multiples of
 the energy of Lebesgue measure on
the unit square.
 (The constants are random, but do not depend on the kernel.)
As an application, we give almost-sure asymptotics for the probability that 
an $\alpha$-stable process approaches within $\epsilon$
of $B[0,1]$, conditional on  $B[0,1]$.

The  upper bound on energy is based on a  strong law
for the approximate self-intersections of the Brownian path.
 
 We  also prove analogous capacity estimates for planar Brownian motion
and for the zero-set of one-dimensional Brownian motion.
\end{abstract}
\vfill
{\bf Keywords:} Brownian motion, capacity, energy,
 occupation measure, local time.
\end{titlepage}
\section{Introduction and main results} \label{sec:main}
\pagestyle{plain}
It is  well-known that for $d \geq 2$, the range of $d$-dimensional Brownian
motion has Hausdorff dimension 2, but its 
2-dimensional measure is almost surely 0. Hausdorff dimension is defined 
via Hausdorff measures, but has an
equally important interpretation (due to Frostman \cite{Fro})
as the critical parameter for
positivity of Riesz capacities.
Exact Hausdorff measure is one much-studied 
means of specifying more precisely  the size of a ``small'' set
(see Taylor \cite{Taylor} for a comprehensive survey in the 
context of random sets);
 exact capacity is a different one, that is directly
relevant to intersections of the small set with other random sets.
 Cieselski and Taylor \cite{CT} found
the exact Hausdorff measure for the trace of Brownian motion in space,
which quantifies to what extent the trace is ``smaller'' than the plane.
 Here we show
 that with probability 1, the spatial Brownian 
  trace has positive capacity exactly in the same kernels as the plane.
Theorem \ref{thm:capequiv.sp} is a quantitative version of this; Theorems 
\ref{thm:capequiv.pl} and \ref{thm:capequiv.Z} give analogous statements for 
planar Brownian motion
and for the zero-set of 1-dimensional Brownian motion, respectively.
 The latter theorem sharpens an integral test due to Kahane and Hawkes.

For a decreasing kernel function
$f : [0,\infty) \rightarrow  [0,\infty]$,
  define the
 {\bf energy} of a Borel measure  $\nu$ on ${\bf R}^d$ with respect to $f$ by
 \[    {\cal E}_f(\nu) =\int_{R^d} \! \int_{R^d }
 f(|x-y|) \, \d\nu(x) \,\d\nu(y)  
 \]
 and the {\bf capacity} of a Borel set
 $\Lambda \subset {\bf  R}^d$ with respect to
 $f$ by
 \[ \Cap_f(\Lambda) = \left[ \inf_{\nu(\Lambda)=1} {\cal E}_f(\nu) \right]^{-1}  . \]
 Thus $ \Cap_f(\Lambda) > 0$ if and only if there exists a Borel measure $\nu$
 supported on $\Lambda$ such that  \newline
${\cal E}_f(\nu) < \infty$. 
  When $f(r) = r^{-\alpha}$, we write  $ \Cap_\alpha $ for
 $ \Cap_f$, and then 
 the ``capacitary dimension" $\sup\{\alpha :  \Cap_\alpha(\Lambda) > 0\}$ 
 of a Borel set $\Lambda$ is equal to its Hausdorff dimension (see, e.g.,
Carleson \cite{Car} or Kahane \cite{Kahane}, page 133).

In the sequel we assume  that all kernel functions
$f$ considered are (weakly) decreasing and satisfy
$\lim_{r \downarrow 0 } f(r) = f(0)$ if this limit is finite.

 Pemantle and Peres \cite{PP}  introduced a notion of
 ``capacity-equivalence", which we specialize to ${\bf R}^d$:
 \begin{dfn} The sets 
$A, B \subset  {\bf R}^d$ are {\bf capacity-equivalent} if there
 exist positive constants $  C_1, C_2  $ such that
 \[ C_1 \Cap_f (B) \leq \Cap_f (A) \leq C_2  \Cap_f (B) 
 \mbox{\rm \hspace{.2in} for all } \, f. \]
 \end{dfn}

Let $(B_t(\omega) \,  : 0 \leq t \leq 1)$ 
be $d$-dimensional Brownian motion started at 0,
and consider its range   
  $B[0,1] = \{x \in {\bf R}^d : B_t
   = x \mbox{ for some } 0 \leq t \leq  1 \}$.
It is classical, and  follows easily from Theorem~\ref{thm:PP}
below (see the discussion around (\ref{eq:cap01})), that
for any kernel $f$, if $m$ denotes  Lebesgue measure on  $[0,1]^2$, then
\[{\cal E}_f(m) \leq C \left[  \Cap_f ([0,1]^2) \right] ^{-1} , \]
where $C$ is an absolute constant.
In particular, $[0,1]^2$ has positive capacity with respect to the kernel 
function $f$ if and only if
\[ \int_{0+}f(r) r\, \d r  < \infty .\]
Theorem~\ref{thm:capequiv.sp} implies, {\em a fortiori},
that with probability 1, the same 
criterion holds for $B[0,1]$ in dimension $d \geq 3$, uniformly over kernels.
\begin{thm} \label{thm:capequiv.sp}
For $d\geq3$, the Brownian trace  
$B[0,1]$ is a.s.\ capacity-equivalent to $[0,1]^2$.
More precisely, with probability~1 there exist random constants 
$C_1, C_2 > 0$ such that
\be \label{eq:capequiv.sp}  
C_1  \Cap_f ([0,1]^2) \leq \Cap_f (B[0,1]) \leq C_2  \Cap_f ([0,1]^2)
\mbox{\rm \hspace{.2in} for all }\, f.
\ee
\end{thm}
In dimension 2, the recurrence of $(B_t)$ leads to a slight modification:
\begin{thm} \label{thm:capequiv.pl}
For any decreasing $f$, denote
 $\widetilde{f}(r) = f(r)\log\frac{1}{r}$. For planar Brownian motion,
  with probability 1
there exist random constants $ C_1, C_2  > 0$ such that 
\be \label{eq:capequiv.pl}
C_1  \Cap_{\widetilde{f}}\,([0,1]^2)
\leq  \Cap_f (B[0,1]) \leq C_2  \Cap_{\widetilde{f}}\,([0,1]^2)
\mbox{\rm \hspace{.2in} for all } \, f. 
\ee
\end{thm}

Our main interest in capacity is that for many stochastic processes, 
particularly
Markov processes (see \cite{Chu}, \cite{FS} and the references therein)
and certain fractal percolation processes
 (see \cite{PP}),
  hitting probabilities of sets are equivalent to their capacities.  

  The next theorem exploits this equivalence,  as well as the  fact
that our almost-sure capacity estimates hold  uniformly over 
  all kernels.  
Aizenman  \cite{Aizenman}
showed that if $[B]$ and $[B^\prime]$ are the traces of
two independent $d$-dimensional
Brownian motions started apart, then
\[
\P\Big[ {\rm dist}([B],[B^\prime]) \, < \epsilon \Big] \asymp 
\left\{
\ba{ll}
\epsilon^{d-4} & \qquad  \mbox{if } d > 4 \\
(\log \frac{1}{\epsilon})^{-1}  
& \qquad \mbox{if } d = 4 
\ea
\right.
\]
as $\epsilon \downarrow 0$.  
(Earlier, 
 Lawler  \cite{Lawler} had obtained precise asymptotics for the analogous
 problem for two random walks on~$\Z^4$. See Albeverio and Zhou \cite{AlbZhou}
 for a recent refinement of Aizenman's estimates.)
Theorem 2.6 of \cite{Pe} contains the following generalization of
Aizenman's result: 
 if $[X^\alpha]$ and $[B]$ denote the traces of an
independent $\alpha$-stable process and Brownian motion,
started apart, then 
\[
\P \Big[ {\rm dist}([B],[X^\alpha]) \, < \epsilon \Big] \asymp
\left\{
\ba{ll}
\epsilon^{d- \alpha - 2} & \qquad  \mbox{if } \alpha < d -  2 \\
(\log \frac{1}{\epsilon})^{-1}  
& \qquad \mbox{if } \alpha  = d - 2 
\ea
\right.  
\]
as $\epsilon \downarrow 0$.

We derive an almost-sure  version of these estimates,
uniform over $\alpha$, conditional on the Brownian motion $B$.
For $0 < \alpha \leq 2$, 
let $\P ^\alpha_x$ be the law of a symmetric $\alpha$-stable process 
$(X^\alpha_t)$ 
in $\R^d$ started at $x$,
 so that
\[ \E ^\alpha_x \, e^{i \lambda \cdot (X_t - x)}
= e^{- |\lambda|^\alpha t} 
\]
for $\lambda \in \R^d$,
and let $f^{(\alpha)}(|x-y|) = c(\alpha)|x-y|^{\alpha -d}$ 
be the corresponding potential density.
We always consider $B$ and $X^\alpha$ to be independent.
Write $[B] = B[0,1]$ and $[X^\alpha] = X^\alpha[0,\infty)$.
\begin{thm} \label{thm:approach}
Suppose $d \geq 3$. 
Let
\bean
m(x, B) & = &  \inf_{y \in [B]}  |x - y| \\
M(x,B) & = & \sup_{y \in [B]}  |x - y|.
\eean
Then for some constants $c_d, c^\prime_d > 0$ the following is true: 
For a.e. Brownian path $B$ and all $x \in \R^d$,
there exists $\epsilon_0 = \epsilon_0(B,x)$
such that, for all $0 < \epsilon < \epsilon_0$,
\[
 c_d \, 
\, M(x, B)^{\alpha - d}
\; \leq  \; 
\frac{\P^\alpha_x \Big[\, {\rm dist}
	\Big([B],[X^\alpha]\Big) < \epsilon \; 
		 \Big|  \, B  \Big]}
	{ \alpha (d - \alpha - 2) \,  \epsilon^{d - \alpha - 2}}
\; \leq \; c_d^\prime  \,  
m(x, B)^{\alpha - d}
\]
for all $0 < \alpha < d - 2$ such that $\alpha \leq 2$, 
and  when $d=3,4$ also
\[
c_d  \,  
\,   M(x, B)^{ - 2}
\;  \leq  \; 
\frac{\P^\alpha_x  \Big[\,{\rm dist}
	\Big([B],[X^\alpha]\Big) < \epsilon \; 
		 \Big|  \, B  \Big]}
	{ (d - 2) \, \left(  \log \frac{1}{\epsilon} \right)^{-1} }
\; \leq \;  c_d^\prime  \,
m(x, B)^{ - 2}
\]
for $\alpha = d - 2$.
\end{thm}

\noindent{\sc Remark:} Note the uniformity in $\alpha$
 in the  statement above. Even for a fixed $\alpha$, the proof of
Theorem \ref{thm:approach}, 
given in section 4, requires estimating the capacity of
a fixed sample path $[B]$ in  infinitely many kernels simultaneously.

Theorems \ref{thm:capequiv.sp} and \ref{thm:capequiv.pl} 
say nothing about which measures supported on $B[0,1]$ 
have low energy with respect to different kernels. It turns out that,
up to a random constant not dependent on the kernel, 
{one measure fits all kernels.}
Let $\mu$  denote the {\bf occupation measure} of $(B_t)$, defined by
\[ \mu(\Lambda)  = \int_0^1  {\bf 1}_{\Lambda}(B_t)\,\d t  \]
for Borel sets $\Lambda \subset {\bf R}^d$.  
Clearly $\mu$ has total mass 1 and is supported
on $B[0,1]$.
Roughly speaking,
for questions of capacity,
$\mu$ plays the same role for $B[0,1]$ that  Lebesgue measure plays
for  $[0,1]^2$.
More precisely, the lower bounds on $ \Cap_f (B[0,1])$
in  (\ref{eq:capequiv.sp}) and
 (\ref{eq:capequiv.pl})   follow directly via
Theorem~\ref{thm:PP} from the next theorem,
 which says that, with probability one, 
the energy of $\mu$ on $B[0,1]$ is bounded by a random constant times the
energy of Lebesgue measure on the unit square, uniformly over kernels. 
\begin{thm} \label{thm:ubint}
With probability one, there exists a $C = C(\omega)$
such that
\be  \label{eq:ubint}
{\cal E}_f(\mu) \leq C\left\{ \ba{lcl}
			 \int_0^1 r f(r)\,\d r  & , & d \geq  3 \\[.05in]
		    	 \int_0^1 r\log \frac{1}{r}f(r)\, \d r  & , & d=2
			  \ea \right.
 \mbox{\rm \hspace{.1in} for all $\, f$.} 
\ee
\end{thm}

A key tool for the proof of the above theorems
is  a simple formula for energy proved in
Benjamini and Peres \cite{BP} (for logarithmic energy) and in
 Pemantle and Peres \cite{PP} (for general kernels), which we state later as
Theorem~\ref{thm:PP}.
As we will show, the upper bounds on  capacity 
 given in Theorems~\ref{thm:capequiv.sp} and \ref{thm:capequiv.pl}
 follow easily from known asymptotics for the volumes of Wiener sausages. 
 The lower bounds on  capacities  are, as we illustrate
 in section~\ref{sec:lb},
 easily proved for {\em fixed}
 kernels, but the fact that, with probability one, these bounds hold  
 uniformly over kernels, is new.   The proofs use
 Theorem~\ref{thm:PP} together with 
  Theorem~\ref{thm:sl}   below.  The proof of 
  Theorem~\ref{thm:approach}, given in section~\ref{sec:approach},
  is similar, and uses the additional 
  deterministic fact 
  that the capacity of an $\epsilon$-sausage  is equivalent to 
  the capacity of the original set with respect to an
  $\epsilon$-smoothed kernel (Proposition~\ref{prop:sauscap}),
  together with the equivalence of capacities and hitting probabilities for 
  stable processes.

For $\sigma > 0$ and $y \in
{\bf R}^d$, define
\[g_\sigma(y) = \exp(-|y|^2/2\sigma^2), \]
 where $|\cdot|$ denotes Euclidean norm.
Let
\bean S_\sigma & = & \int_0^1 \d t_1 \int_0^1 \d t_2 \;
g_\sigma(B_{t_1}-B_{t_2}) \\
& = & \int_{R^d} \! \int_{R^d} g_{\sigma}(x-y) \,\d \mu(x)\, \d \mu(y) \, .
\eean
When suitably scaled, $S_\sigma$ may be interpreted as  measuring the
``approximate self-intersections"
of the Brownian path.  The case $d=2$ of the following theorem follows from
  Varadhan's renormalization of $S_\sigma$
  (see section \ref{sec:last}
for details).
\begin{thm}[A strong law for approximate self-intersections] \label{thm:sl} \,
For $d \geq 3$,
\[
\frac{S_\sigma}{\sigma^2} \rightarrow \frac{4}{d-2}
\mbox{\hspace{.2in} as $\sigma \downarrow 0$ , a.s.}  \]
In dimension 2,
\[
\frac{S_\sigma}{\sigma^2 \log\frac{1}{\sigma}} \rightarrow 4
\mbox{\hspace{.2in} as $\sigma \downarrow 0$ , a.s.}  \]
\end{thm}

To explain the connection to the energy estimates in Theorem \ref{thm:ubint},
we start with the observation that the ratio $\mu(Q)/{\mbox{\rm side}(Q)^2}$
cannot be uniformly bounded as $Q$ ranges over all cubes,
 since the 2-dimensional
Hausdorff measure of the Brownian trace vanishes. 
The \linebreak $\mu$-weighted average of this ratio, taken
over the collection $\dn$ of all dyadic cubes $Q$ of side $2^{-n}$, is
\be \label{eq:sumsum} 
4^n \sum_{Q \in \dn }  \mu(Q)^2 .
\ee
Theorem \ref{thm:sl} implies
that, in dimension $d \geq 3$, these weighted averages are bounded uniformly
in $n$. (See the inequality (\ref{eq:sumsq}) in Subsection \ref{sec:ubintpf}.)
Theorem \ref{thm:PP} is then used to express the energy of $\mu$
as a positive linear combination of the  averages (\ref{eq:sumsum}),
 and thus to compare it to the energy of 
Lebesgue measure on the unit square. 
\subsection{ The zero set}
We have  analogous results for the zero set
of one-dimensional Brownian motion,
 \linebreak $Z = \{t \in [0,1] : B_t = 0 \}$. 
These results are technically
easier than the corresponding ones for the Brownian trace,
 and led us to the latter.
It is classical that $Z$ a.s.\ has 
Hausdorff dimension 1/2 (again, with zero measure in that dimension),
 so here a natural 
comparison set is the ``middle-1/2 Cantor set''
\[ K = \left\{ \sum_{n=1}^\infty b_n 4^{-n} : b_n = 0,3 \right\} .\]
  $K$ is a standard example of a set of Hausdorff dimension $1/2$, that has
positive and finite measure in that dimension.
\begin{thm} \label{thm:capequiv.Z} The Brownian zero-set $Z$ is a.s.\
capacity-equivalent to the middle-1/2 Cantor set $K$.
More precisely, with probability one 
there  exist random $C_1, C_2 > 0$, such that
\be \label{eq:capequiv.Z}
C_1  \Cap_f(K) \leq \Cap_f (Z) \leq C_2  \Cap_f(K)
\mbox{\rm \hspace{.2in} for all } \, f.
\ee
\end{thm}
Let $(\ell(t) : 0 \leq t \leq 1)$ be Brownian local time at zero, normalized so
that, by results of  L\'{e}vy, $\ell(t) $ has the same law as the  
running maximum
$\max_{\tau \leq t}B_\tau$. We abuse notation slightly and also let
$\ell$ denote the 
measure, supported on $Z$, for which it is the distribution function.

 The  lower bound on $\Cap_f (Z)$ in 
(\ref{eq:capequiv.Z}) is implied by the following energy estimate:
\begin{thm} \label{thm:ubint.Z}
With probability one there exists $C = C(\omega)$ such that:
\be \label{eq:ubint.Z}
{\cal E}_f(\ell) \leq C \int_0^1f(r)r^{-1/2}\, \d r ,
\mbox{\rm \hspace{.1in}  for all } \, f .
\ee
\end{thm}

In the first (1968) edition of \cite{Kahane}, Kahane 
 established that, for a fixed  $f$ of ``positive type'', 
  finiteness of the integral 
 in (\ref{eq:ubint.Z}) is 
 sufficient for the Brownian zero set $Z$ to a.s.\ have 
 positive capacity with respect to $f$ (see 
\cite{Kahane} page 236, Theorem 2). This is the first ``exact capacity''
result we are aware of. Hawkes (\cite{Hawkes} Theorem~5)
 proved the converse (finiteness of the integral is necessary for positive
 capacity) under a slightly stronger assumption (log-convexity) on the 
 kernel $f$. In view of   the expression (\ref{eq:capK})
for the capacity of $K$, Theorem \ref{thm:capequiv.Z}
 is a uniform version of this result of Kahane and Hawkes;
 it also shows that the side conditions on the kernel are not needed.
In the last section we describe a different random set that
illustrates why  the uniformity in the
kernel is not automatic.

\section{Upper bounds on capacities }
\label{sec:ubcap}

The following representation of energy from \cite{PP} is basic for most of
the results in this paper. Its proof is based on a trick from \cite{BP}.
Let $\dn$  denote the collection of
all dyadic cubes \newline  $Q =
[j_1 2^{-n},(j_1+1)2^{-n}) \times \ldots \times [j_d2^{-n},(j_d+1)2^{-n})$ for
$ (j_1,\ldots,j_d) \in {\bf Z}^d$.
\begin{thm} [\cite{PP}, Theorem~3.1]  \label{thm:PP}Let
$f : [0,\infty) \rightarrow  [0,\infty]$ be a weakly decreasing function.
Then for any Borel measure $\nu$ supported  on
the unit cube $[0,1]^d$,
\be \label{eq:sumbyp}
 {\cal E}_f(\nu) \asymp
 \sum_{n=0}^\infty (f(2^{-n}) - f(2^{1-n}))\sum_{Q \in \dn} \nu(Q)^2 ,
\ee
where  $\asymp$ means that the ratio of the two quantities
 is bounded between
two positive constants depending only on $d$.
\end{thm}
\noindent{\sc Remark:} The proof of this in \cite{PP} assumes that
$f(0+)=\infty$ and that $\nu$ has no atoms,
 but these assumptions can be avoided as
long as $f(0)=f(0+)$. If $\nu$ has atoms at the points $\{x_j\}_{j \geq 1}$,
 then there is a contribution of $\sum_j f(0) \nu(x_j)^2$ to the energy 
${\cal E}_f(\nu)$  coming from the diagonal.
 On the right-hand side of~(\ref{eq:sumbyp}), we get the same contribution.

We first note an easy general 
upper bound on capacity, which is essentially the same as
Theorem IV.2 in Carleson  \cite{Car}.
Let $N_n(\Lambda)$ 
be the number of dyadic cubes  $Q \in \dn$ (as defined in 
Theorem~\ref{thm:PP}) that intersect a Borel set $\Lambda \subset {\bf R}^d$.
  Then there is a constant $c>0$, depending only on the ambient
dimension $d$, such that for any 
probability measure $\nu$ supported on $\Lambda$, and any kernel $f$, we have 
\bea 
{\cal E}_f(\nu) & \geq &
c\sum_n \left(f(2^{-n}) - f(2^{1-n})\right)\sum_{Q \in \dn} \nu(Q)^2 \nonumber \\
 & \geq & c \sum_n  \left(f(2^{-n}) - f(2^{1-n})\right) N_n(\Lambda)^{-1}. \nonumber
   \eea 
Therefore
\be  \label{eq:capub} \Cap_f(\Lambda) \leq 
c^{-1} \left[  \sum_n \left(f(2^{-n}) - f(2^{1-n}) \right) 
N_n(\Lambda)^{-1} \right]^{-1} .
\ee

If for some $c$, the set 
$\Lambda$ carries a positive measure $\nu$ such that 
$\nu(Q) \leq c N_n(\Lambda)^{-1}$ for all $Q \in \dn$ and all $n$, then 
this bound is sharp (up to  a constant factor independent of $f$). 
Thus we get
\be  \label{eq:cap01} 
\Cap_f ([0,1]^2) \, \asymp \,
	  \left[ \sum_n  (f(2^{-n}) - f(2^{1-n})) 4^{-n}\right]^{-1}
\, \asymp \, \left[ \int_0^1 f(r) r \, \d r \right]^{-1}
\ee
and similarly, for the middle-half Cantor set
\be \label{eq:capK} 
\Cap_f (K) \, \asymp \, 
  \left[ \sum_n  (f(2^{-n}) - f(2^{1-n})) 2^{-n/2}\right]^{-1} 
 \, \asymp \,  \left[ \int_0^1 f(r) r^{-1/2} \, \d r \right]^{-1} ,
 \ee
 where  $\asymp$ means that the ratio of the two sides is bounded above and
    below by  positive absolute constants.  The minimum energies 
   are  attained within a constant factor 
  by  Lebesgue measure in the case 
   of $[0,1]^2$, and, for $K$, by the measure that makes
  the digits  $(b_n)$~$\sim$~i.i.d~Bernoulli(1/2), 
  when $K$ is represented as 
$ \left\{ \sum_{n=1}^\infty b_n 4^{-n} : b_n = 0,3 \right\}$.

{\bf Proof of Theorems~\ref{thm:capequiv.sp} and \ref{thm:capequiv.pl}
- upper bound: }
Strong laws for volumes of Wiener  sausages (see \cite{LG92} Chapter VI and the
references therein) imply that, with probability one, there 
exist random $C_1, C_2 \in 
(0, \infty)$ such that for all $n$, 
\be \label{eq:saus}
\ba{ccccccl}
C_1 & \leq & {N_n ( B[0,1]) \over   4^n } & \leq C_2 & \mbox{for} & d\geq 3 \\
[.15in]
 C_1 & \leq & {n \cdot N_n(B[0,1])\over 4^n } & \leq C_2 & \mbox{for} & d = 2 .
   \ea \ee
Substituting the above into (\ref{eq:capub}) and comparing 
with (\ref{eq:cap01}) gives, with probability one, 
\[  \Cap_f (B[0,1])  \leq C(\omega) 
		\left\{ \ba{l}
		     \Cap_{\widetilde{f}}\,([0,1]^2) \;\; , \;\; d=2 \\
		      \Cap_f ([0,1]^2) \;\; , \;\; d\geq3
		      \ea \right.
\mbox{\rm \hspace{.2in} for all }\, f, \]
where $\widetilde{f}$ is  defined in the statement of 
Theorem~\ref{thm:capequiv.pl}.   $\Cox$

{\bf Proof of Theorem~\ref{thm:capequiv.Z} - upper bound: } 
We need an analog of 
(\ref{eq:saus}) for $Z$. This is provided by 
Kingman's \cite{Kingman} construction of local time, which we 
sketch here for the Brownian case.  Recall L\'{e}vy's classical result   
(see, e.g., \cite{RY} page 447)
\be \label{eq:Lev}
\delta^{1/2}\widetilde{N}_\delta \rightarrow 
\left(\frac{2}{\pi}\right)^{1/2}\ell(1) 
\mbox{\hspace{.2in}, as $\delta \downarrow 0$, a.s.}
\ee
where $\widetilde{N}_\delta$ is the number of maximal intervals $I_j$ of
$[0,1]\backslash Z$ having length 
greater than $\delta$. Now if 
\[ Z^\delta = \{ u \in [0, 1] : B_t=0 
\mbox{ for some $t$ with } |u-t|<\delta/2 \}\; ,
\]
and $m$ denotes Lebesgue measure on ${\bf R}^+$, then, using the fact 
that $m(Z)=0$  a.s., we obtain
$$
m(Z^\delta)  =  \sum_j \left[m(I_j) \wedge \delta\right] + O(\delta) \, ,
$$
where the sum extends over all maximal intervals in $[0,1] \setminus Z$.
By Fubini's theorem,  this sum can be written as 
$ \int_0^\delta  \widetilde{N}_\epsilon \, \d \epsilon  + O(\delta)$.
Together with 
(\ref{eq:Lev}), this  implies that 
\[  \delta^{-1/2} m(Z^\delta) \rightarrow 
2\left(\frac{2}{\pi}\right)^{1/2} \ell(1) 
\mbox{\hspace{.2in} as $\delta \downarrow 0$, almost surely.}
\]
 Thus for suitable absolute constants  $c_1, c_2 > 0$,
there almost surely exists 
 a random integer $n^*$, such that
\be \label {eq:NZ}
c_1 \ell(1) \leq 2^{-n/2} N_n(Z) \leq c_2 \,  \ell(1) \mbox{ \hspace{.2in}
 for all}  \, n \geq n^* \, .  
\ee
The upper bound on $\Cap_f(Z)$ now follows from the general upper bound
 (\ref{eq:capub})
and the estimate (\ref{eq:capK}). 
$\Cox$

\section {Lower bounds on capacities} \label{sec:lb}
 
 {\bf Remark: } For a {\em fixed} kernel, it is easy to see  that finiteness of
 the integral on the right hand side of  (\ref{eq:ubint}) or 
 (\ref{eq:ubint.Z}) implies that the left hand side is finite.
 We  show this for  (\ref{eq:ubint}) in the case $d=3$; the other proofs
 are similar.
  Recall that for any non-negative Borel function $h : \R ^3 \rightarrow \R$,
\be \label{eq:green}
 \E \int_0^\infty h(B_t)\,  \d t = \frac{1}{2\pi}
		\int_{R^3}  h(x)\frac{\d x}{|x|}\, ,
\ee
  where $|\cdot|$ is the  Euclidean norm and  $\d x$ denotes Lebesgue measure.
 By the Markov property,  we have
  \bean
  \E{\cal E}_f(\mu)
  & \leq & 2\E\int_0^1 f(|B_t|) \, \d t  
\eean
Since $f$ is monotone decreasing, 
$\, f(|x|)  \leq  f(|x|)  \one_{\{|x| \leq 1\}} + f(1)$.
Invoking (\ref{eq:green}), we get
$$
 \E{\cal E}_f(\mu) \, \leq \, 
 2 \left(  \frac{1}{2\pi}
  \int_0^1 f(r) \cdot 4\pi r^2  \, \frac{\d r}{r} 
		+ f(1) \right) \\
 \, \leq \,  (4+4) \,  \int_0^1 f(r) \, r \, \d r   ,
$$
  where  the last step used the monotonicity of $f$ again.
$\Cox$

\subsection{The Brownian trace} \label{sec:ubintpf}
{\sc Proof of Theorem \ref{thm:ubint}: } Recall that $\dn$ is the collection
of dyadic squares of side $2^{-n}$.
For $\sigma = 2^{-n}$ we have, by  the definition
of $S_\sigma$, that 
$$
S_{\sigma} \geq \sum_{Q \in \dn} 
    \int_Q \int_Q g_{\sigma} (x-y) \,  \d\mu(x) \, \d \mu (y) \, .
$$ 
All the integrands on the right hand side are bounded below by a positive 
constant $c = c(d)$ which does not depend on $n$.
 Hence by Theorem~\ref{thm:sl},
 there is a random constant 
$C^\prime= C^\prime (\omega) $ such that, with probability one, for all~$n$ 
\be \label{eq:sumsq}
\sum_{Q \in \dn} \mu(Q)^2 
\, \leq \, c^{-1} S_{2^{-n}}
\, \leq \,  C^\prime  \left\{ \ba{lcl}
		4^{-n} & , & d \geq 3 \\
		n4^{-n} & , & d=2
		\ea \right. .
\ee
Thus, by Theorem \ref{thm:PP}, with probability 1
\bean
{\cal E}_f(\mu) 
& \leq &   C^\prime c_1  \sum_{n=n_0}^\infty
	\left( f(2^{-n}) - f(2^{1-n}) \right)
	\left\{ \ba{lcl}
	4^{-n} & , & d \geq 3 \\
	n4^{-n} & , & d=2
	\ea \right. , 
\eean
where $n_0 = n_0(\omega)$ is defined by
 $ 2^{-n_0} \geq  \mbox{diameter }(B[0,1]) > 2^{-n_0-1}$, and $c_1$
 depends only on $d$. 
Since $f$ is monotone decreasing, by adjusting $C^\prime$ we may 
replace $n_0$ by 1 in the above sum,
and we obtain (\ref{eq:ubint}) after a summation by parts. $\Cox$

\subsection{The zero set}
We first prove a proposition, which, loosely speaking, will play the role
that  Theorem~\ref{thm:sl} did in the previous proof.
Recall that $\ell(\cdot)$ denotes local time at 0.
\begin{prop}
 \label{prop:ltsl}  
Consider the quadratic variation of $\ell$ at scale $\delta$:
\bean
L_\delta
   & = & \sum_{j = 0}^{\lceil \delta^{-1} \rceil}
      \left[\ell((j+1)\delta) -
       \ell(j\delta)\right]^2 .
\eean
With probability~1, there exists a random 
$C = C(\omega)$ such that
\be \label{eq:qv}  
      L_\delta \leq C \delta^{1/2} 
\mbox{\rm \hspace{.2in} for all  }\, \delta > 0 .  
\ee
\end{prop}

{\sc Proof: } We consider separately the summands for odd and even $j$ in 
$L_\delta$. Denote one-dimensional Brownian motion by $B_t$.
For fixed $\delta>0$,  
let $j_1, j_2, \ldots, $ be a left-to-right enumeration
of all the \underbar{odd} $j \geq 1$
such that  $B_t=0$ for some $t$ in the interval $[(j-1)\delta \, , \, j \delta]$.
Let $ M(\delta) := \max \{i \, : \,  j_i \delta \leq 1+ \delta \} $  
be the number of these intervals which intersect $[0,1]$.

Define stopping times
 $T_i = \inf\{t\in [(j_i-1)\delta \, ,  \, j_i \delta] \, : \, B_t = 0\}$,
and let $X_i:= \ell (T_i + \delta) -\ell (T_i).$
The strong Markov property at the times $T_i$ implies that, for fixed $\delta$,
the variables
$\{X_i\}_{i \geq 1}$
 are i.i.d.\   with the law of $\ell(\delta)$, which is the same
as the law of $|B_{\delta}|$. In particular $X_i^2$ have mean $\delta$ and 
exponentially decaying tails. Thus the partial sums
$Y_k(\delta) := \sum_{i=1}^k X_i$ satisfy
\be \label{eq:ld}
\P \Big(Y_k(\delta) > 2k \delta \Big) \leq e^{-ck}
  \; \, \mbox{for some constant} \; c> 0 .
\ee
By the argument leading to
(\ref{eq:NZ}), with probability 1 there exists a $\delta^* = \delta^*(\omega)$
such that
\be \label{eq:ltsl.1}
  M(\delta) \leq c_2 \ell(1) \delta^{-1/2}
\mbox{\rm \hspace{.1in}  for all } \, \delta < \delta^* ,
\ee
with $c_2$ an absolute constant.
 
Denote $Y^{(n)} := Y_{M(2^{-n})} (2^{-n})$.
Since $k=k(n) = c_2 \ell (1) 2^{n/2}$ is eventually larger than $n$, we see that
\bean
\P \Big[ Y^{(n)} \, > \, 2 c_2 \ell (1) 2^{-n/2} \; \mbox{\bf i.o.} \,\Big]
& \leq &  \, \\
\P \Big[ M(2^{-n}) > c_2 \ell (1) 2^{n/2}\;  \mbox{\bf i.o.}\, \Big] & + &
\P \Big[\, \mbox{for infinitely many} \; n, \; \exists k >n \,:  \, 
   Y_k(2^{-n}) > 2k 2^{-n}\Big] \, .
\eean
The first probability in the sum vanishes by (\ref{eq:ltsl.1}), and the second by
(\ref{eq:ld}) and Borel-Cantelli.
Thus a.s.\ there is a random constant $A=A(\omega)$
such that $Y^{(n)} \leq A 2^{-n/2}$ for all $n$.
Now $Y^{(n)}$ is an upper bound for the sum over all odd indices $j$ in the 
quadratic variation $L_{2^{-n}}$, and the even indices are handled similarly.
Consequently $2^{n/2} L_{2^{-n}} $ is a.s.\ bounded by a random constant. 

To go from the powers of $1/2$ to general $\delta$, observe that
any interval $I$ can be covered by three shorter dyadic intervals,
say $J_1, J_2, J_3$. Clearly 
$\ell(I)^2 \, \leq  \, 3 (\ell(J_1)^2+\ell(J_2)^2+\ell(J_3)^2)$.
Therefore, if $2^{1-n} > \delta \geq 2^{-n}$ then
$L_{\delta}  \leq 6 L_{2^{-n}}$. This concludes the proof.
  $\Cox$

{\sc  Proof of Theorem \ref{thm:ubint.Z}: }
Follow the proof of Theorem~\ref{thm:ubint} given in 
section~\ref{sec:ubintpf},   
replacing $\mu$ by $\ell$  and using 
Proposition~\ref{prop:ltsl}.
 $\Cox$

\section{ Probabilities of $\epsilon$-approach }
\label{sec:approach}
In this section we prove Theorem~\ref{thm:approach}.
The next deterministic proposition
 states that the capacity of an $\epsilon$-sausage  is equivalent to
     the capacity of the original set with respect to an
       $\epsilon$-smoothed kernel.
More precisely,
given a kernel function $f$ and $\epsilon > 0$,
let  
\[
\overline{f}(\epsilon) = 
\epsilon^{-d} d  \int_0^{\epsilon } f(s) \,s^{d-1} \,
\d s , \]
and define 
\[ f_\epsilon (r) =
\left\{
\ba{ll}
f(r) & \qquad \mbox{ if } r \geq \epsilon   \\
\overline{f}(\epsilon) &  \qquad \mbox{ if } r <  \epsilon 
\ea \right. .
\]
Note that  $f_\epsilon$ is decreasing, since $f$ is.
Also, $\overline{f}(\epsilon) < \infty$ provided that
$\Cap_f (\R^d) >0$, which we may always assume.

For a Borel set
$\Lambda \subset \R^d$, we denote the $\epsilon$-sausage about $\Lambda$ by
\[ \Lambda_\epsilon = \{x : |x-y|  < \epsilon
\; \; \mbox{ for some } y \in \Lambda\}.
\]
Recall  that ``$\asymp$''
(``is comparable to'') means that the  two quantities
are within finite positive constant multiples of each other,
the constants depending only on the dimension~$d$.
Similarly, the expression ``$a \lequiv b$" will
mean ``$a \leq c_d \, b$".
We also use the notation $Q \in \dn $ for dyadic cubes introduced at
the beginning of section~\ref{sec:ubcap}.
\begin{prop} \label{prop:sauscap}
For any Borel set $\Lambda \subset \R^d$, kernel function $f$,
and $\epsilon > 0$, we have
\labeq{eq:sauscap}
 \Cap_f (\Lambda_\epsilon) \asymp \Cap_{f_\epsilon}(\Lambda).
  \ee
  \end{prop}
  {\sc Proof: }
It clearly suffices to prove the proposition for compact $\Lambda$.
We first show that the left-hand side of (\ref{eq:sauscap}) is,
up to a constant factor,
greater than the right.
  Given a probability measure $\nu$ on $\Lambda$,
 it is natural to smooth it by convolving with normalized Lebesgue measure
 on a ball of radius $\epsilon$. It will be even easier to control a discrete
 version of this convolution. Choose $m_\epsilon $ and  $n_\epsilon $ so that
\bean
 &2^{-m_\epsilon }\:< \: \epsilon
  \: \leq \:  2^{-m_\epsilon +1} \; & \mbox{ and } \\
  &\sqrt{d} 2^{-n_\epsilon } \: < \: \epsilon
  \: \leq \: \sqrt{d} 2^{-n_\epsilon +1} \, . &
\eean 
Observe that the definition of $\overline{f}(\epsilon)$ and the monotinicity
of $f$ imply that
\be \label{eq:fbar}
 \overline{f}(\epsilon) \asymp 
 \sum_{n \geq m_\epsilon}  f(2^{-n})  2^{d(m_\epsilon-n)} \, .
\ee 
  
Define a smoothed probability measure
  $\nu_\epsilon$ by
  \[ \d \nu_\epsilon \,\Big|_Q
  = 2^{n_\epsilon d} \, \nu(Q) \, \d x \,\Big|_Q  ,
  \qquad \mbox{ for } Q \in {\cal D}_{n_\epsilon}, \]
  where $\d x$ denotes Lebesgue measure.

  Suppose $\nu$ is supported on $\Lambda$;
 then $\nu_\epsilon$ is supported on $\Lambda_\epsilon$.
 Note that for every $n$ we have 
$$
\sum_{Q \in \dn} \nu_\epsilon (Q)^2 \leq 2^{d(n_\epsilon-n)}
   \sum_{Q \in {\cal D}_{n_\epsilon}}
    \nu (Q)^2 \, ;
$$
 indeed for $n  \geq n_\epsilon$ the two sides are clearly equal,
while for $n  < n_\epsilon$ the inequality follows from Cauchy-Schwarz,
since every $Q \in \dn$ is the union of $2^{d(n_\epsilon-n)}$ cubes in
  ${\cal D}_{n_\epsilon}$.
  Thus using  (\ref{eq:sumbyp}) to expand
  ${\cal E}_f (\nu_\epsilon)$ gives
  \labea{ea:sauscap}
  {\cal E}_f ( \nu_\epsilon )
  & \lequiv &
\sum_{n<m_\epsilon}
	   \Big(f(2^{-n}) - f(2^{1-n}) \Big)\sum_{Q \in \dn} \nu (Q)^2 \\[1ex]
	      & &+  \sum_{n \geq m_\epsilon} \Big( f(2^{-n}) - f(2^{1-n}) \Big)
    2^{d(n_\epsilon-n)}
   \sum_{Q \in {\cal D}_{n_\epsilon}}
    \nu (Q)^2 \, .\nonumber 
\eea
 Since $2^{d n_\epsilon} \asymp 2^{d m_\epsilon}$, 
 by (\ref{eq:fbar}) the last line is comparable to
$$
  \Big(\overline{f}(\epsilon) -  f(2^{1-m_\epsilon}) \Big)
  \sum_{Q \in {\cal D}_{n_\epsilon}}
    \nu (Q)^2 \, .  
$$
Invoking (\ref{eq:sumbyp}) again, we infer that
\be \label{eq:compcaps}
{\cal E}_f (\nu_\epsilon) \lequiv  {\cal E}_{f_\epsilon} (\nu) \,.  
\ee
(The reverse inequality $\gequiv$ also holds, but we will not need it.)
  The  asserted inequality 
$\Cap_f (\Lambda_\epsilon)^{-1} \lequiv \Cap_{f_\epsilon}(\Lambda)^{-1}$
now follows by taking the infimum in (\ref{eq:compcaps}) 
 as $\nu$ ranges over probability measures on $\Lambda$.

To obtain the reverse inequality, we use a Borel-measurable mapping
 $\pi: \Lambda_\epsilon \rightarrow \Lambda$,
 which moves every point by at most
$\epsilon$. For instance, $\pi(x)$ can be defined as the lexicographically
minimal $y \in \Lambda$ such that $|y-x| \leq \epsilon$.
 
Suppose  that $\nu$ is a probability measure on
$\Lambda_\epsilon$, and consider the projected measure 
 $\nu \pi^{-1}$ on $\Lambda$.
As before, we have
\bea
  {\cal E}_{f_\epsilon} (\nu  \pi^{-1})
     & \asymp &
      \sum_{n<m_\epsilon}
	\Big(f(2^{-n}) - f(2^{1-n})\Big) \sum_{Q \in \dn} \nu  \pi^{-1} (Q)^2
	  \nonumber \\[1.8ex]
	     & &  \; \: + \;     
	   \Big(\overline{f}(\epsilon) -  f(2^{1-m_\epsilon}) \Big)
  \sum_{Q \in {\cal D}_{m_\epsilon}}
    \nu \pi^{-1}(Q)^2 \, . \label{ea:scpi} 
\eea
Now for each cube ${Q \in \dn}$, the
preimage $\pi^{-1}(Q)$ is contained in the union of the cubes
 $Q^\prime \in D_n$ such that $\dist(Q^\prime , Q ) < \epsilon$.
If $n \leq m_\epsilon$, then there are at most $5^d$ such cubes $Q^\prime$,
and hence by Cauchy-Schwarz,
$$
 \nu  \pi^{-1} (Q)^2   \leq 
  5^d \sum_{Q^\prime \in \dn}
                 \nu (Q^\prime\, \cap \, \pi^{-1} Q)  ^2 \, . 
$$
Therefore for $n \leq m_\epsilon$,
\be  \label{ea:scpi2}
\sum_{Q \in \dn} \nu  \pi^{-1} (Q)^2  \, \leq \,
    5^d \sum_{Q' \in \dn}
       \Big(  \sum_{Q \in \dn}
	  \nu (Q^\prime \cap \, \pi^{-1} Q)^2 \Big) \nonumber \\
  \, \leq \,  5^{d} \sum_{ Q^\prime \in  {\cal D}_n}
     \nu(Q^\prime)^2 \, .
\ee
Combining (\ref{ea:scpi}) and (\ref{ea:scpi2}), we get
\bea
  {\cal E}_{f_\epsilon} (\nu  \pi^{-1}) & \lequiv &
      \sum_{n<m_\epsilon}
	\Big(f(2^{-n}) - f(2^{1-n})\Big)\sum_{Q \in \dn} \nu   (Q)^2
	  \nonumber \\[1.8ex]
	     & &  \; \: + \;     
	   \Big(\overline{f}(\epsilon) -  f(2^{1-m_\epsilon}) \Big)
  \sum_{Q \in {\cal D}_{m_\epsilon}}
    \nu (Q)^2 \, . \label{ea:scpi3} 
\eea
 On the other hand, we can use Cauchy-Schwarz to bound the energy
${\cal E}_f  (\nu)$ from below:
\bea     
    {\cal E}_f  (\nu) &\gequiv&
     \sum_{n<m_\epsilon}
	\Big(f(2^{-n}) - f(2^{1-n})\Big) \sum_{Q \in \dn} \nu   (Q)^2
	  \nonumber \\[1.8ex]
	     & &  \; \: + \;     
	\sum_{n \geq m_\epsilon} \Big(f(2^{-n}) - f(2^{1-n})\Big) 2^{m_\epsilon-n}
  \sum_{Q \in {\cal D}_{m_\epsilon}} 
    \nu (Q)^2 \, . \label{ea:scpi4} 
\eea
				  
   By using (\ref{eq:fbar}) to compare (\ref{ea:scpi3})
  and  (\ref{ea:scpi4}), we see that
    \[
      {\cal E}_{f_\epsilon} (\nu  \pi^{-1})
       \lequiv
	 {\cal E}_f (\nu),
	 \]
	 and taking the infimum over probability measures $\nu$
  on $\Lambda_\epsilon$ completes the proof.
 $\Cox$

Next, we recall the  well-known quantitative version
of the classical equivalence between the capacity of a set  
and its probability of being hit by a  stable process.
As in the introduction, let $\P^\alpha_x$ denote the law of 
 a symmetric $\alpha$-stable
process $(X^\alpha_t)$ started at $x \in \R^d$ with
potential density
$f^{(\alpha)}(|x-y|) =  c(\alpha) \,  |x-y|^{\alpha - d}$
and trace $[X^\alpha]$.
\begin{prop}(see, e.g., \cite{Tay67} Lemma 2, or 
   \cite{Pe} Proposition~3.2) \label{prop:hitprob}
Let $\Lambda$ be any Borel subset of $\R^d$, and
suppose there are positive numbers $k$ and $K$ such that
$k \leq f^{(\alpha)}(|x-y|) \leq K$ for all $y \in \Lambda$.  Then
\[
k \, \Cap_{f^{(\alpha)}} (\Lambda)
\; \leq \P^\alpha_x 
	\Big[ [X^\alpha] \cap \Lambda \neq \emptyset \Big]
\; \leq \; K  \, \Cap_{f^{(\alpha)}} (\Lambda).
\]
\end{prop}

{\sc Proof of Theorem~\ref{thm:approach}:}
Recall the notation $\overline{f}(\epsilon)$ and $f_\epsilon$
introduced at the beginning of this section.
The proof begins similarly to that of Theorem~\ref{thm:ubint}.
By Theorem~\ref{thm:sl},  for some 
fixed  constants $c$ and  $c^\prime > 0$,
with probability 1 there exists $n_* = n_*(\omega)$
such that
\labeq{eq:app1}
 c  4^{-n} \; \leq \; S_{2^{-n}} \; \leq \; c^\prime  4^{-n}
\qquad \mbox{for } n > n_*.
\ee
By (\ref{eq:sumbyp}) and (\ref{eq:sumsq}),
we have
\be \label{eq:twosums}
{\cal E}_{f^{(\alpha)}_\epsilon}(\mu)
\lequiv
\left( \sum_{n \leq n_*} \; + \:  \sum_{n > n_*} \right)
\left(f^{(\alpha)}_\epsilon(2^{-n}) - f^{(\alpha)}_\epsilon(2^{-n+1}) \right)
\, S_{2^{-n}}.
\ee
Assume that   $\epsilon < 2^{-n_*}$.
Then the first sum is clearly $\leq f^{(\alpha)}(2^{-n_*})$. 
On the other hand, a simple integration shows that
\labeq{eq:app2}
\overline{f^{(\alpha)}}(\epsilon) \; = \;  \frac{c(\alpha)}{\alpha} \,
 \epsilon^{\alpha - d} \; \mbox{ for } \,  \epsilon>0.
\ee
Assume now that $\alpha < d - 2$.
Substituting (\ref{eq:app1}) into the second sum in (\ref{eq:twosums}),
summing by parts  (as in the proof of Theorem~\ref{thm:ubint}),
and letting $\epsilon \downarrow 0$ 
shows that 
\[
{\cal E}_{f^{(\alpha)}_\epsilon}(\mu)
\lequiv  \frac{c(\alpha)}{\alpha (d - \alpha - 2)}
\, \epsilon^{2+\alpha - d}
\]
for all $\epsilon$ less than some $\epsilon_0(\omega)$.
So, by Proposition~\ref{prop:sauscap},
\[ \Cap_{f^{(\alpha)}} \left( [B]_\epsilon \right)
\; \asymp
\;  \Cap_{f^{(\alpha)}_\epsilon} \left( [B] \right)
\;  \gequiv
\;   \frac{\alpha (d - \alpha - 2)}{c(\alpha)}
\, \epsilon^{d - 2 - \alpha }
\]
for $\epsilon < \epsilon_0$.
Since ${\rm dist}([X^\alpha],[B]) < \epsilon$ if and only if
$X^\alpha$ hits $[B]_\epsilon$, the above estimate
and Proposition~\ref{prop:hitprob}
	establish the desired lower bound 
	on
	$\P^\alpha_x\Big[\,{\rm dist}([X^\alpha],[B]) < \epsilon 
	\; \Big|\, B \Big]$.
A similar calculation handles the case $\alpha = d - 2$.
	 The proof of the upper bound is entirely analogous, using the
  general upper bound on capacity (\ref{eq:capub}) and the
	 strong law for volumes of Wiener sausages
	 alluded to  above (\ref{eq:saus}) instead of Theorem~\ref{thm:sl}.
	 $\Cox$

	\section{ Proof of the strong law for $S_\sigma \;$
		  (Theorem \protect\ref{thm:sl}) }
	We prove  Theorem~\ref{thm:sl} only for the case $d\geq3$.
      Our elementary  
	method also works with only minor 
	modifications for $d=2$, but  since this case  follows from Varadhan's
	renormalization, which has received 
	at least four proofs (\cite{Var, Rosen86a, LG85, Yor86}),
	we omit it here.
	{\bf Throughout this section, 
	we assume $d\geq 3$.}

	The argument follows classical lines: Estimate the first two moments,
	use Chebyshev's inequality to obtain convergence along a subsequence, and 
	interpolate. However, showing that the variance of $S_{\sigma}$ is of lower 
	order than the squared mean requires some care, so we include the details.
	\subsection{Moment estimates}
	Define the joint probability  densities \newline 
	$p(t_1,\ldots,t_k; x_1,\ldots,x_k)$ by
	\[ \P(B_{t_1}\in A_1, \ldots,B_{t_k}
	\in A_k) = \int_{A_1 \times \ldots \times A_k} \d x_1 \ldots \d x_k \,
	p(t_1,\ldots,t_k;x_1,\ldots,x_k) , \]
	$A_i \subset {\bf R}^d$ Borel.
	\begin{prop}
	\label{prop:es1}
	\be \label{eq:es1}
	 \E S_\sigma  =  
		\frac{4}{d-2}\sigma^2 \;  + \;  \Theta_d(\sigma) 
	\ee
	where
	\[
	 \Theta_d(\sigma)  =  \left\{ \ba{lcl}
				O(\sigma^3) & , & d=3 \\
				O(\sigma^4\log\frac{1}{\sigma})  & , & d=4\\
				O(\sigma^4)     & , & d\geq 5
				\ea \right. 
	\] 
	as $\sigma \downarrow 0$.
	\end{prop}
	{\sc Proof:} By definition,
	\bean \E S_\sigma & = &
	2 \int_{0\leq t_1 \leq t_2 \leq 1} \d t_1 \,\d t_2 \,\int_{(R^d)^2}
	 \d x_1 \, \d x_2 \,
	p(t_1,t_2;x_1,x_2) \, g_\sigma(x_1-x_2) \\
	& = & 2 \int_0^1 \d s \, \frac{1-s }{(2\pi s)^{d/2}} \int_{R^d} \d y \,
	\exp \left[-\frac{1}{2}|y|^2(\frac{1}{s}+\frac{1}{\sigma^2}) \right] \, ,
	\eean
	after changing variables $s\equiv t_2-t_1$ and $y\equiv x_2-x_1$ and integrating
	out first $ x_1$ and then $t_1$. Therefore
	\bean
	\E S_\sigma & = & 2 \int_0^1 \d s \, (1-s) 
	\left(\frac{\sigma^2}{\sigma^2 + s}\right)^{d/2}\\
	& = & 2 \sigma^d\int_0^1  \frac{\d s}{(\sigma^2 + s)^{d/2}}
	-  2 \sigma^d\int_0^1  \frac{s \, \d s}{(\sigma^2 + s)^{d/2}}.
	\eean
	One readily checks that the first term equals the right-hand side of
	(\ref{eq:es1}), while the second
	is easily bounded using
	\bean
	\int_0^1  \frac{s\, { \d s} }{(\sigma^2 + s)^{d/2}}
	 & \leq & \int_0^1  \frac{\mbox{\boldmath $d$}s}{(\sigma^2 + s)^{\frac{d}{2}-1}}
	 \mbox{\hspace{.2in}} \Cox
	  \eean

	 \begin{prop} [The second moment] \label{prop:es2}
	 \be \label{eq:es2}
	 \E S_\sigma^2 = \left( \frac{4}{d-2} \sigma^2 \right)^2
				   \; + \; \Theta_d^\prime(\sigma)  
	\ee
	where
	\[
\Theta_d^\prime(\sigma)  =  \left\{ \ba{lcl}
						O(\sigma^5) & , & d=3 \\
					 O(\sigma^6\log\frac{1}{\sigma})  & , & d=4\\
					 O(\sigma^6)     & , & d\geq 5
						\ea \right. 
	\]
	as $\sigma \downarrow 0$.
	\end{prop}

	In the calculations below, we always have $s , s_i, t \geq 0$.
	We will repeatedly use the following bound:

	\bea 
	\int \! \! \int_{s+t \leq 1} \d s \,\d t  \, 
	 \frac{s}{(\sigma^2 + s)^{d/2}} \frac{1}{(\sigma^2 +t)^{d/2}}
	& \leq & \int_0^1   \frac{ \d s}{(\sigma^2 + s)^{d/2-1}}
		 \int_0^1   \frac{\d t }{(\sigma^2 + t)^{d/2}} \nonumber \\
	& = & \left\{ \ba{lcl}
			O(\sigma^{-1})	        & , & d=3  \\[.03in]
			O(\sigma^{-2} \log\frac{1}{\sigma})  & , & d=4 \\[.03in]
			O(\sigma^{6 - 2d})	   & , & d \geq 5
			\ea  \label{eq:a1}
		\right.
	\eea
	Call these orders of magnitude $\Psi_d(\sigma)$.  Note that
	\be \label{eq:a2} 
	\sigma^{2d} \Psi_d(\sigma) = \Theta^\prime_d(\sigma).
	\ee

	{\bf Proof of Proposition \ref{prop:es2}:}
	\bean
	 \E S_\sigma^2 & = & 8 
	\int \!\!\! \int  \!\!\! \int \!\!\! \int_{0\leq t_1\leq  \ldots  \leq t_4 \leq 1}
	  \d t_1 \ldots \d t_4  
	\int \!\!\! \int  \!\!\! \int \!\!\! \int_{(R^d)^4} \d x_1 \ldots \d x_4 \,
	   \, p(t_1,\ldots,t_4 ; x_1,\ldots,x_4) \\[.1in]
	   & & \times \left\{g_\sigma(x_1-x_2)g_\sigma(x_3-x_4)
	   + g_\sigma(x_1-x_3)g_\sigma(x_2-x_4)
	    + g_\sigma(x_1-x_4)g_\sigma(x_2-x_3)\right\} \\[.15in]
	  & = & 8 ( I_1 + I_2 + I_3) ,
	   \eean
	 say.
	 The calculations below
	 show that 8$I_1$ is equal to the right side of
	(\ref{eq:es2}), and that the other integrals are of the smaller order.
	 The latter fact makes intuitive sense:  as $\sigma\downarrow 0$, the major
	 contribution to each $I_i$ comes from the region of the time simplex  where
	 the  path increments being 
	 weighted by $g_\sigma$ have small time increments.
	  But for $I_2, I_3$, this requires that at least {\em three\/} time-increments
	be small simultaneously, putting us in a corner of the  simplex and so
	losing  powers of $\sigma$ asymptotically.  

	{\bf Estimation of $I_1$: }
	Changing variables $s_i \equiv t_{i+1}-t_i$ and 
	$y_i \equiv x_{i+1}-x_i$,  and integrating 
	out two unweighted space-time  increments,
	\bea \label{eq:I1} 
	I_1 & = & 
	\int \!\! \! \int_{s_1+s_3\leq 1}\d s_1 \,\d s_3 \, 
	\frac{(1-s_1-s_3)^2}{2} \nonumber\\[.1in]
	  &  & \times  \frac{1}{(2\pi s_1)^{d/2}}
	\int_{R^d}\d y_1 \,\exp \left[-\frac{1}{2} |y_1|^2
	 \left(\frac{1}{s_1} + \frac{1}{\sigma^2}\right)\right] \nonumber \\[.1in]
	   &  & \times   \frac{1}{ (2\pi s_3)^{d/2}}
	    \int_{R^d}\d y_3 \,\exp \left[-\frac{1}{2} |y_3|^2
	      \left(\frac{1}{s_3} + \frac{1}{\sigma^2}\right)\right] \nonumber\\[.15in]
	 & = & \sigma^{2d} 
	\int \!\!\! \int_{s_1+s_3\leq 1}\d s_1 \,\d s_3 \,\frac{(1-s_1-s_3)^2}{2}
	\frac{1}{(\sigma^2 + s_1)^{d/2}} \frac{1}{(\sigma^2 + s_3)^{d/2}}.
	 \eea

	 Expanding (\ref{eq:I1}) and using (\ref{eq:a1}) and  (\ref{eq:a2}),
	\bea 
	I_1 
	& = & \frac{1}{2} \sigma^{2d}
		 \int \!\! \! \int_{s_1+s_3\leq 1}\d s_1 \, \d s_3 \,
		 \frac{1}{(\sigma^2 + s_1)^{d/2}} \frac{1}{(\sigma^2 + s_3)^{d/2}}
		\; + \; \Theta^\prime_d(\sigma). \label{eq:expandI1}
	\eea
	 To handle the
	first term, 
	\bea
	 &&\int_0^1 \d s_1 \, \frac{1}{(\sigma^2 + s_1)^{d/2}}
	 \int_0^{1-s_1} \d s_3 \,  \frac{1}{(\sigma^2 + s_3)^{d/2}} \nonumber \\[.13in]
	& = & \int_0^1 \d s_1 \, \frac{1}{(\sigma^2 + s_1)^{d/2}}
	\frac{2}{d-2} \left( \sigma^{2 - d} 
	- (\sigma^2 + 1 - s_1)^{1 - \frac{d}{2}} \right)  . \label{eq:I1spat}
	\eea
	The first term of this is
	\[
	 \left( \frac{2}{d-2} \right)^2 \sigma^{4-2d} 
	+ O(\sigma^{2-d})\]
	while the absolute value of the second (negative) term in (\ref{eq:I1spat})
	is bounded 
	 by integrating on $[0,1/2]$ and $[1/2,1]$ separately:
	$$
	\int_0^1 \d s_1 \, \frac{1}{(\sigma^2 + s_1)^{d/2}} 
	\frac{1}{(\sigma^2 + 1 - s_1)^{\frac{d}{2}-1}}
	 \leq  \Psi_d(\sigma) \,
	$$
	with room to spare.
	Multiplying everything by $8\cdot \frac{1}{2} \sigma^{2d}$ gives the
	right-hand side of (\ref{eq:es2}). 

	 {\bf  Estimation of $I_2$: } With the same change of variables 
	 $s_i \equiv t_{i+1} - t_i$ and  $y_i \equiv x_{i+1} - x_i$, we integrate
	 out $y_0$ and $s_0$ to obtain
	\bean
	I_2 & = & 
	\int \!\! \!\int \!\! \!\int_{s_1+s_2+s_3\leq 1}
	\d s_1 \, \d s_2 \, \d s_3 \, (1-s_1-s_2-s_3)\\[.1in]
	& & \times  \frac{1}{(2\pi s_1)^{d/2}}
	\frac{1}{(2\pi s_2)^{d/2}}
	\int \!\! \!\int_{(R^d)^2}\d y_1\, \d y_2  \, 
	\exp -\frac{1}{2}\left(\frac{|y_1|^2}{s_1}
	  + \frac{|y_2|^2}{s_2} + \frac{|y_1+y_2|^2}{\sigma^2}\right)\\[.1in]
	 & & \times      \frac{1}{(2\pi s_3)^{d/2}}
	\int_{R^d} \d y_3 \,
	    \exp -\frac{1}{2} \left(
	   \frac{|y_3|^2}{s_3} + \frac{|y_2+y_3|^2}{\sigma^2}\right)  .
		  \eean
	 Changing variables $t \equiv s_1+s_2 , z \equiv y_1+y_2$ and integrating
	 out $y_1$ and $s_1$, we  get
	 \bea
	  I_2 & = &\int \!\! \! \int_{t+s_3\leq 1}\d t\,  \d s_3 \, (1-t-s_3)t 
		\nonumber \\[.1in]
	  & & \times  \frac{1}{(2\pi t)^{d/2}}
	   \int_{R^d}\d z \, \exp \left[-\frac{1}{2} |z|^2
	     \left(\frac{1}{t} + \frac{1}{\sigma^2}\right) \right] 
		\nonumber \\[.1in]
	      & & \times    \frac{1}{(2\pi s_3)^{d/2}}
	      \int_{R^d}\d y_3 \,
	   \exp -\frac{1}{2} \left(\frac{|y_3|^2}{s_3}
	     +  \frac{|y_2+y_3|^2}{\sigma^2}\right) . \label{lastline}
	\eea
	We bound the last factor above
	(line~(\ref{lastline}))  
by noticing that it is
\bean
& = & (2 \pi \sigma^2)^{d/2}
\, \int_{R^d} p(s_3;y_3) \, p(\sigma^2; -y_2-y_3) \, \d y_3  \\
&  = &  (2 \pi \sigma^2)^{d/2}
\, p(s_3 + \sigma^2; -y_2)  \\
& \leq &  \frac{\sigma^d}{(\sigma^2 + s_3)^{d/2}}.
\eean
	Thus
	\bean
	I_2 & \leq & 
	      \sigma^{2d} \int \!\! \!\int_{t+s_3\leq 1}\d t \, \d s_3 \, (1-t-s_3)t
	     \frac{1}{(\sigma^2 + t)^{d/2}} \frac{1}{(\sigma^2 +s_3)^{d/2}} \\[.13in]
	     & \leq &  \sigma^{2d} 
	\int \!\! \! \int_{t+s_3\leq 1}\d t \, \d s_3 \,
	     \frac{1}{(\sigma^2 + t)^{\frac{d}{2}-1}} \frac{1}{(\sigma^2 +s_3)^{d/2}} 
	\\[.13in]
	& = & \Theta^\prime_d(\sigma) 
	\eean
	by  (\ref{eq:a1}) and (\ref{eq:a2}).

	{\bf Estimation of $I_3$: } Similarly,
	\bean
	I_3 & = & \int \!\! \!\int \!\! \!\int_{s_1+s_2+s_3\leq 1}
	\d s_1 \, \d s_2 \, \d s_3 \,(1-s_1-s_2-s_3)
	\\[.1in]
	& & \times  \frac{1}{(2\pi s_2)^{d/2}}
	\int_{R^d}\d y_2 \, \exp \left[-\frac{1}{2} |y_2|^2
	  \left(\frac{1}{s_2} + \frac{1}{\sigma^2}\right) \right]\\[.1in]
	  & & \times    \frac{1}{(2\pi s_1)^{d/2}}  \frac{1}{(2\pi s_3)^{d/2}}
	\int_{(R^d)^2}\d y_1 \, \d y_3 \,
	     \exp -\frac{1}{2} \left(\frac{|y_1|^2}{s_1}
		  + \frac{|y_3|^2}{s_3} + \frac{|y_1+y_2+y_3|^2}{\sigma^2}\right)  .
	  \eean
	 Changing variables $t \equiv s_1+s_3 , z \equiv y_1+y_3$ and integrating
	 out $y_1$ and $s_1$, we have
	 \bean
	 I_3 & = & \int \!\! \! \int_{t+s_2\leq 1}\d t \, \d s_2 \, (1-t-s_2)t  \\[.1in]
	  & & \times  \frac{1}{(2\pi s_2)^{d/2}}
	\int_{R^d}\d y_2 \, \exp \left[-\frac{1}{2} |y_2|^2 \left(\frac{1}{s_2} + 
	\frac{1}{\sigma^2}\right) \right]\\[.1in]
	& & \times    \frac{1}{(2\pi t)^{d/2}}
	   \int_{R^d}\d z \,
	 \exp -\frac{1}{2} \left(\frac{|z|^2}{t}
	   +  \frac{|z +y_2|^2}{\sigma^2}\right) .
	\eean
As at (\ref{lastline}), the last factor above is 
bounded by $ \frac{\sigma^d}{(\sigma^2 +t)^{d/2}}$,
	and we obtain
	\bean
	  I_3 & \leq &  \sigma^{2d}  \int \!\! \!
		\int_{t+s_2\leq 1}\, \d t \, \d s_2 (1-t-s_2)t
	  \frac{1}{(\sigma^2 + s_2)^{d/2}} \frac{1}{(\sigma^2 +t)^{d/2}} \\[.13in]
	  & \leq &  \sigma^{2d} \int\!\! \!
		\int_{t+s_2\leq 1}\, \d t \, \d s_2
	    \frac{1}{(\sigma^2 + s_2)^{d/2}} \frac{1}{(\sigma^2 +t)^{\frac{d}{2}-1}} 
	\\[.13in]
	& = & \Theta^\prime_d(\sigma) 
	\eean
	 by (\ref{eq:a1}) and (\ref{eq:a2}).

	\subsection{Almost-sure convergence}
	We need the following deterministic lemma.
For any  Borel measure  $\nu$ on ${\bf R}^d$ and $\sigma>0$, define
\[ S_\sigma(\nu) = \int_{R^d} \!\int_{R^d }
g_\sigma (|x - y|) \,  d\nu(x) \, d\nu(y).
\]
Thus $S_\sigma = S_\sigma(\mu)$.
\begin{lem}	\label{lem:mon}
For any Borel measure $\nu$ on  ${\bf R}^d$, the quantity  
$\sigma^{-d}S_\sigma(\nu)$ is monotone decreasing in $\sigma$.
In particular,  $\sigma^{-d}S_\sigma$ is a.s.\
monotone decreasing in $\sigma$.
\end{lem}
{\sc Proof:}
Let $\, \widehat{ } \,$ denote the Fourier transform, so that for $\xi \in {\bf R}^d$
\bean
\widehat{g_\sigma}(\xi) & = & 
(2\pi)^{-d/2}\int_{{\bf R}^d} g_\sigma(x)e^{-i\xi \cdot x}\,\d x \\
	& = & \sigma^d e^{-\sigma^2 |\xi|^2/2}  .
\eean
Then, by Plancherel's formula,
\bean
\sigma^{-d}S_\sigma(\nu) &  = & 
\sigma^{-d}\int_{{\bf R}^d}  \widehat{g_\sigma}(\xi) |\widehat{\nu}(\xi)|^2 
\,d \xi \\ 
& = &  \int_{{\bf R}^d}  e^{-\sigma^2 |\xi|^2/2}  
|\widehat{\nu}(\xi)|^2 \,d \xi .
\eean
and the lemma clearly follows. $\Cox$

{\bf Proof of Theorem \ref{thm:sl} :}

 By Propositions~\ref{prop:es1} and \ref{prop:es2},
\[ {\bf Var }(S_\sigma) = \left\{
		 \ba{lcl}
		 O(\sigma^5) & , & d=3 \\
		  O(\sigma^6\log \frac{1}{\sigma}) & , & d=4 \\
		   O(\sigma^6) & , & d\geq 5
		   \ea \right. .
\]
By Chebyshev's inequality, for any $\epsilon > 0$,
$$
\P\Big(|S_\sigma -\E S_\sigma|>\epsilon \sigma^2 \Big) 
  \leq   \epsilon^{-2} \sigma^{-4} {\bf Var }(S_\sigma) = O(\sigma).
$$
The right hand side is summable as $\sigma$ runs over the sequence  
$\sigma_n = n^{-2}$,
so by Borel Cantelli and  Proposition~\ref{prop:es1}, 
\be
\label{eq:subseq1}
{\sigma_n}^{-2} S_{\sigma_n} \rightarrow {4 \over d-2}
 \mbox{\rm \hspace{.2in}as $n  \rightarrow \infty $, a.s.}
\ee
Now  for arbitrary  positive $\sigma < 1$, choose $n$ such that 
$\sigma_{n+1} < \sigma \leq \sigma_n$. Then by Lemma \ref{lem:mon},
\[ \ba{ccccl}
\sigma_{n}^{-d}\,S_{\sigma_{n}} & \leq  &\sigma^{-d}\,S_\sigma 
&\leq & \sigma_{n+1}^{-d}\,S_{\sigma_{n+1}}
\ea
\]
so that
\[    \ba{ccccl}
(\sigma /{ \sigma_{n}})^{d-2}\,  
\sigma_{n}^{-2} \, S_{\sigma_{n}} & \leq  &\sigma^{-2}\,S_\sigma 
&\leq & (\sigma /{ \sigma_{n+1}})^{d-2}\, \sigma_{n+1}^{-2}\,S_{\sigma_{n+1}}\, .
\ea
\]

Thus $\sigma^{-2}\,S_\sigma $ is sandwiched between two
expressions which tend to ${4 \over d-2}$ as $\sigma \downarrow 0$,
and we're done. 
$\Cox$

\section{Concluding remarks} \label{sec:last}

\begin{enumerate}

\item  
The following example shows that the
uniformity in capacity-equivalence statements for random sets
is not automatic.
Consider the
random Cantor set $\Lambda$ in $[0,1]$ constructed as follows.
For each $k \geq 1$, pick a random integer $n_k$ uniformly
in the interval $[3^k+k, 3^{k+1}-k]$, with all picks independent;
define $\Lambda$ to be the set of all sums
$\sum_{n=1}^\infty a_n 4^{-n}$ with 
$$a_n= 
     \left\{ \ba{lcl}
       			0 & \mbox{\rm for} & n \in (n_k-k,n_k]\\[.03in]
      			0,1,2,3 & \mbox{\rm for} &n \in (n_k,n_k+k] \\[.03in]
			      0,3 & \, &  \mbox{otherwise.} 
  			\ea \right. 
$$
Then it is not hard to check that for fixed $f$, with probability one,
$\Cap_{f}(\Lambda) > 0$ if and only if 
 $ \int_0^1f(r)r^{-1/2}\, \d r < \infty$. 
(See \cite{Pe2} for details.)
However,
$\Lambda$ is {\em not} capacity-equivalent
to the middle-half Cantor set; indeed there   exists 
 a random kernel $f^* $ (depending on the sample $\Lambda$)
that satisfies this integrability
condition  but gives
 $\Cap_{f^*}(\Lambda) = 0$.

\item In 1969 Varadhan \cite{Var} proved that, in dimension two, 
$\sigma^{-2}(S_\sigma - \E  S_\sigma)$ converges a.s. to a well-defined 
random variable.  This clearly implies the planar case of
our Theorem~\ref{thm:sl}. 
Varadhan's renormalization has received many proofs and extensions. (See
\cite{Var, Rosen86a, LG85, Yor86, Dyn85, Dyn88} as well as 
Chapter VIII of \cite{LG92} and 
the bibliographical notes there.)  Rosen \cite{Rosen86b},
Remarks II-III, gives detailed 
calculations which  are close in spirit to ours 
(though more Fourier-analytic),
and which could probably be extended to $d \geq 3$ to yield our 
Theorem \ref{thm:sl}.
While self-intersection local time exists in dimension 3 as well as 
dimension 2, there seems 
 to be no analogue of Varadhan's almost-sure renormalization there.  
However, Yor \cite{Yor85} shows that the renormalized $S_\sigma$ converge 
{\em in law}
 for $d=3$.   Yor establishes that 
 $\sigma^{-3} (\log \frac{1}{\sigma})^{-1/2}(S_\sigma - \E S_\sigma)$
 converges in law (to a Gaussian) as $\sigma \downarrow 0$;
 this seems tighter than the estimate $\Var   [S_\sigma] = O(\sigma^5)$
given in (\ref{eq:es2}).
\item  \it Does Brownian motion in three-space almost 
surely have the property that
{\em all} of the orthogonal projections to planes of its trace are 
capacity-equivalent to each other? \rm \newline
\item Let $B$ and $B^\prime$ be two independent standard Brownian
motions in $\R^3$. The ``fractal percolation"
methods of \cite{PP} and \cite{Pe}, 
which are based on the results of Lyons \cite{Ly}, imply the following:
for any fixed kernel $f$, the capacity of the
intersection $\Cap_f(B[0,1] \cap B^\prime[0,1])$ 
is almost surely positive if $\Cap_f([0,1]) >0$;
otherwise with probability 1  the intersection has capacity 0  in this kernel.
However, these methods do not indicate if this holds uniformly in the kernel.
\begin{quote}
\it Is the intersection $B[0,1] \cap B^\prime[0,1]$ of two independent 
 Brownian traces in $\R^3$, almost-surely capacity-equivalent to $[0,1]$?
\rm
\end{quote}
\end{enumerate}

\noindent{\bf Acknowledgement:} We are indebted to Russell Lyons, who
first alerted us to the importance of the order of quantifiers in capacity
estimates for random sets.

\renewcommand{\baselinestretch}{1.0}

\end{document}